\documentclass{commatDV}

\usepackage{DLde}

\title{Arkady Onishchik: on his life and work on supersymmetry}

\author[Dimitry Leites]{Dimitry Leites${}^{a, b}$}

\affiliation{${}^a$New York University Abu Dhabi, 
Division of Science and Mathematics, P.O. Box 129188, United Arab
Emirates; 
${}^{b}$Depart\-ment of Mathematics, Stockholm University, SE-106 91
Stockholm, Sweden
\email{dl146@nyu.edu; mleites@math.su.se}%
}

\date{}



\abstract{Selected stories about the life of A.\,L.~Onishchik, and a~ review of his contribution to the classification of non-split supermanifolds, in particular, supercurves a.k.a. superstrings; his editorial and educational work. A~brief overview of his and his students' results in supersymmetry, and their impact on other researchers. 

Several open problems growing out of Onishchik's research are presented, some of them are related with odd parameters of deformations and non-holonomic structures (non-integrable distributions) on supermanifolds important in physical models, such as Minkowski superspaces and certain superstrings.
}

\VOLUME{30}
\NUMBER{3}
\firstpage{1}
\DOI{https://doi.org/10.46298/cm.9337}

\begin{paper}

This is a~disjoint union of my recollections of certain aspects of the life of one of my teachers with glimpses of ``local color'' of that time.

\textbf{Arkady Lvovich Onishchik (14.11.1933--12.02.2019)} was quiet-looking, but actually rather passionate and charismatic man. I want to at least try to describe these features of Onishchik, and his contribution to the study of supercurves or (for physicists) superstrings and other supermanifolds of prime interest (homogeneous superspaces whose underlying is an Hermitian symmetric space). For his other mathematical achievements, see \cite{AJ}, \cite{AOb}.

\textbf{Hobbies}. To be able to read inscriptions on coins he collected, he learned --- to the extent needed for this --- Arabic, Chinese, Hebrew, Japanese, Urdu and a~number of European languages. From his childhood, he liked to study languages and he learned a~good deal more than numismatics required: he was able to speak and write in several of the above-listed languages. 

A.L. used to tell his students interested in languages (e.g., E.~Vishnyakova who spoke several) that to seriously study mathematics one has to know at least 5 languages (English, French, German, Italian and Russian).

Onishchik's widow (I mean: his third wife) did not agree with my claim above, but was unable to precisely formulate her version of how many languages Onishchik mastered fluently (``perhaps, six... judging by the experience in the countries we visited''). She thought I was exaggerating saying he was able to speak Chinese, not only read. However, being acquainted with rudiments of Mandarin (which I mostly forgot, regrettably), I can testify that he \textit{was} able to speak Chinese, at least a~bit. In particular, A.L. told me why in his translation into Russian of the book by Chern Shiing-Shen, whose family name nowadays many, after the Americans, pronounce as if ``Chern'' is written in English, and whose name nobody (except for Chinese-speaking people) dares to pronounce at all, 
he wrote (in Cyrillic letters) ``Chzhen' Shen'-Shen": this was according to the accepted rules of transliteration of Chinese words (correctly written Ch\'en X\v{i}ngsh\=en in pinyin encoding, i.e., in Latin characters).

\textbf{On Onishchik's charisma and ability to ``influence people'' (as D.~Carnegie would say)}. When the daughter of his second wife was ca 14 years old, she used to attend Math Circle for pupils at the Moscow State University. Once, helping organizers of the Circle, I saw her at some lecture for pupils, and (not having seen her before) identified her by her last name when somebody addressed her. By now, I forgot her and her mother's Spanish last names which I remembered at that time. In the Soviet Union one sometimes --- for various reasons --- took his/her mother's last name, so I assumed the daughter was common with Onishchik. 

The girl looked intelligent and solved some problems on the spot, so the next day, when I met A.L., I told him, being sure he will be pleased to hear anything positive about his nice daughter: ``I met your daughter at a~ math lecture for kids, she performed very well!'' 

A.L. asked, puzzled: ``What daughter?'' (I did not know at that time that he had a~much older daughter by his first wife.) When I told A.L. her first and  last names, he retorted: ``This is my wife's daughter!'' 

``But she looks your splitting image: the posture and movements!'' 

``No!!! Her eyes are brown, and mine are grey!'' 

``But her posture and way of speaking are exactly yours!'' (I recalled the ethologists' term ``imprinting''.)

\textbf{Taking a~risk}. Once, on the writing desk in his flat I saw the book by I.\,R.~Shafarevich ``Socialism as a~phenomenon of world history'' (in Russian) recently published  by a~notoriously anti-Soviet publisher (YMCA-Press, 1977). 
The book was interesting, especially being forbidden, although certain interpretations in it are wrong or doubtful. Onishchik gave it to a~curious me to read at home, without any comments then or later, at the time he was already on the authorities's ``black list'' for having signed ``the letter by 99'' briefly mentioned without any explanation in \cite{AOb}. (I highly recommend Neretin's  breath-catching investigation (currently only in Russian, regrettably) of the history and circumstances related to this letter, and the list of 99 names, mostly of distinguished mathematicians but naive people, lured into this affair by skilled provocators, see \protect{\url{https://www.mat.univie.ac.at/~neretin/misc/1968plus.html}.}) 

By that time we were acquainted for several years, Onishchik actively learned from me whatever I knew of supersymmetry, but  such trust in his freshly graduated student was still unusual for a~ supervisor. 

Let me explain to those who did not have close experience with certain realities of the Soviet Union of that time (like me, although I lived there), and who are unable to appreciate  this trust in full. Onishchik was risking his job and, in the worst --- unlikely, but possible --- case his freedom. Such extreme punishments for reading a~prohibited book were applied very seldom, and mostly in the province, but not only in the cases where reading was aggravated by other ``crimes'', like signing a~``wrong" letter. 

Later, I knew people punished in this way. One instance happened in Petrozavodsk with my friend M.~Serov. He was invited for an \lq\lq interview" by a~ KGB officer who asked Serov if he knew A.~Lavut (the editor of the underground \lq\lq Chronicle of current events"). Serov answered: \lq\lq Sure, he is my close friend\rq\rq.  Serov was never involved in any dissidents' activity himself and was left free after this \lq\lq interview", but was forced to retire with a~half of otherwise expected pension.  With at least 3 dependents, this was a~heavy blow; he died soon after having thus tacitly added his name to those variously punished openly, by a~ court ruling, e.g., sentenced to hard labor camps and exile, see \url{https://en.wikipedia.org/wiki/Chronicle_of_Current_Events}.

For the other person I knew (one of Onishchik's Ph.D. students, V.~Golitsyn), to lose --- under pressure by  KGB ``curators'' of the University --- his Lecturer position in Kalinin (now Tver) was a~ totally unexpected blessing. Being thus expelled, he applied to a~ factory to work as a~ programmer. He barely dared to own up to the head of HR what for he was previously fired: politics. The reaction was unexpected: ``Here, half of the workers did a~term, some of them for murder. Do not talk nonsense and stop worrying, just go to work." He told me that his salary then was higher than at the University, and there were almost no ``voluntary, but compulsory" meetings or yearly ``potato pickings''. (Recall that scientific researchers, students, and low-rank professors had to spend several days or weeks in collective farms harvesting or performing some other dirty job, mainly useless (the results were destined to rot, mostly) except for demonstrating loyalty to the regime.) Soon Golitsyn became an active dissident, and eventually was fired from the factory as well  (under same type of pressure as before, and against the will of his immediate superiors who were happy with his performance at work). Then, his colleagues from the factory, and their wives, were not afraid to help his family --- unlike his colleagues at the University wherefrom he was fired before...

\textbf{How Onishchik got interested in ``super''}.   F.\,A.~Berezin was the first of my scientific advisors, I knew him from school. He was a~Sr.~ Scientific Researcher at the Chair of the Theory of Functions and Functional Analysis of the Department of Mechanics and Mathematics (Mekh-Mat) of Moscow State University. He decided (mistakenly, as it turned out) that he will securely shelter me from his negative --- in the eyes of Mekh-Mat bosses --- influence, if he will not have accepted me as his official student after my second year of study --- when every sophomore had to find a~scientific advisor and be affiliated with a~ Chair (of Algebra, Logic, etc.) depending on the future specialization. (About Berezin's life at that time, read the book \cite{Shif},  its expounded version \cite{KMLT}, and the papers \cite{KNV}, \cite{N}.)

Having thus explained to me that he did not dare to spoil my career, Berezin immediately made a~phone call to Onishchik  and asked him to be my official scientific advisor.  Onishchik taught the course of algebra to a~ stream of ca 100 listeners consisting of several groups, to one of which I belonged, so he at least had seen my face.

During the meeting that followed, Onishchik asked me what was I doing under Berezin's guidance. I answered that my task, as I understood it, was to define an analog of Grothen\-dieck's scheme for what in 1971 was called the graded-commutative $\Zee/2$-graded ring or algebra (which nowadays, after the 1974 talk by Wess and Zumino (\cite{WZ}) is known as supercommutative super ring or superalgebra) and this task was related to physics, more precisely to a~uniform description of bosons (particles with integer spin), and fermions (particles with half-integer spin) --- the only two types of particles considered  in the (mainstream) models of the modern quantum physics. It would be the affirmative answer to Berezin's conjecture 
\begin{equation}\label{BerConj}
\begin{minipage}[l]{12cm}
``There should exist an analog of Calculus in which the elements of the Grassmann algebra play the role of functions", 
\end{minipage}
\end{equation}
namely: the superscheme and its smooth or analytic version would be the --- elusive at that time --- ``something'' on which the elements of the Grassmann algebra play the role of functions. 

Onishchik said that the problem looked interesting, ``but since it is completely unclear when will you be able  to give the  required definition, whereas I have to grade your course work in early spring next year, here is another problem for you''. The extremely polite Arkady Lvovich did not say ``\textbf{IF} you will be able'', and formulated a~problem which he could solve in several minutes: ``describe spherical functions on some manifold of low dimension''.

The task Onishchik suggested looked boring to me and completely devoid of challenge, unlike Berezin's task through which I perceived overwhelming opportunities and perspectives that would open as soon as the desired definition will have been obtained. I gave the definition of what now is called  a~\textit{superscheme} at the beginning of my third year at the University. 

Next, to have Calculus, one had to define a smooth version of the superscheme. 
Today, it is unclear why instead of saying ``do the same with underlying manifold instead of a scheme'' Berezin and I discussed this, to and fro, for almost 2 years after I have given the definition.
For circumstances that psychologically prepared me to be able to foresee what was unclear to many researchers  despite Berezin's public talks in the late 1960s and through the 1970s, see Preface to the book~\cite{Lsos}. 

For reasons completely unclear now, experts in algebraic geometry whom Berezin tried to get interested in the task, refused to think on his problem. Although the spectrum of these superalgebras --- an analog of manifold on which the elements of the Grassmann algebra play the role of functions --- was indeed the same as the spectrum of the quotient of the superalgebra modulo the ideal generated by the anticommuting indeterminates (so, geometrically,  we \textit{ostensibly} gain nothing) these experts overlooked larger group of symmetries of the superscheme as compared with the symmetry group of the underlying scheme of the quotient algebra, to say nothing about the supergroups of symmetries of the superscheme. 

Moreover, when I solved Berezin's problem --- gave the definition of the superscheme ---   Berezin himself did not accept it: he could not believe that the \textit{spectrum} of the Grassmann algebra $G(n)$ with $n\geq 1$ generators (the set on which $G(n)$ is the algebra of functions)  is just one point.  
Berezin conceded to accept my definition only in August 1974, several months after the Wess-Zumino talk~ \cite{WZ} and after my definition of the superscheme had been published \cite{L0}. 

Therefore, I am forever grateful to A.\,L.~Onishchik who in November of 1972 said ``I do not see what Berezin did not like in your definition. It is correct. Write it down as a~short note and I'll submit it to ``Uspekhi" (Matematicheskih Nauk) as a~Communication of Moscow Mathematical Society.'' That is how the first definition of superscheme (\cite{L0}) was published. The Communications of Moscow Mathematical Society, a~part of the volume of ``Uspekhi" started to be translated cover-to-cover in ``Russian Mathematical Surveys'', only from 1975 on. This encouraged certain people to try to claim priority. Attempts to ``improve'' the definition  of supermanifold or superscheme, or just to rewrite it without due reference but with mistakes, continue till now; thanks to irresponsible or unqualified referees  such papers and even books are being published. For detailed arguments explaining what is wrong in some of (rather numerous) ``improved definitions'', see criticism of ``alternative'' definitions of supermanifolds in  \cite[\S~4.8]{Mo}. Regrettably, these arguments did not stop vacuous or wrong papers and even books on super topic with the same mistakes as those discussed in \cite[\S~4.8]{Mo} to keep appearing later on.

At the end of 1974, my definition of the superscheme was reformulated in a~ joint paper with Berezin to define the smooth supermanifold, but not in terms of sheaves (as was natural for superschemes and supervarieties, especially if the ground field is finite), but in the language of charts and atlases. Berezin said this formulation would be more understandable for physicists, whom he perceived at that time as the main, if not only, target audience. (These approaches are equivalent if both can be applied; the one with charts and atlases is more adequate in the infinite-dimensional case, even in the absence of ``super", see \cite{Mo}, while one can not use it for algebraic manifolds or supervarieties over finite fields, when the only possible languages are the one of sheaves and ringed spaces, or that of the functor of points a.k.a. families.)

A bit after my note \cite{L0} had been published, Wess and Zumino delivered a~talk \cite{WZ} that highly impressed physicists --- and justly so. In this talk they introduced the word (not the precise term yet) \textbf{supersymmetry} and --- much more explicitly than in several  previous works where Lie superalgebras appear with more or less conscious applications to theoretical physics --- showed the usefulness of supersymmetry in field theory. Some historians of science count up to a~half dozen pioneer works, but their authors do not understand the importance and meaning of their own works even now, see Prefaces in \cite{DSB}. The importance of these works published earlier than \cite{WZ}, and containing seedlings of what is now called \textbf{supersymmetries}, was (over)appreciated only after the talk \cite{WZ}. Today's understanding of various applications of supersymmetry had to wait for several decades. 

I do not know if Wess and Zumino were aware that \lq\lq From somewhere in the 1950s on, John Wheeler  repeatedly urged people who were interested in the quantum-gravity program  to understand the structure of a~mathematical object that he called \textit{Superspace} \cite{Wh}. The intended meaning of \textit{Superspace} was that of a~set, denoted by $S(\Sigma)$, whose points faithfully correspond to all possible Riemannian geometries on a~given three-manifold $\Sigma$\rq\rq, see \cite{Giu}. Thus, Wheeler (known not only as Feynman's teacher, but also for having coined the terms \lq\lq black hole" and \lq\lq wormhole") used the term \lq\lq superspace" in a~sense absolutely different from what is now customary to use by everybody, except geometro-dynamists: 

 \begin{quotation}\lq\lq The stage on which the space of the Universe moves is certainly not space itself. Nobody can be a~stage for himself; he has to have a~larger arena in which to move. The arena in which space does its changing is not even the space-time of Einstein, for space-time is the history of space changing with time. The arena must be a~larger object: \textit{Superspace} \dots It is not endowed with three or four dimensions --- it is endowed with an infinite number of dimensions.\rq\rq (J.A. Wheeler: \textit{Superspace}, Harper's Magazine, July 1974, p. 9.)\end{quotation}

It is regrettable that the same term has two completely different meanings because both meanings can meet in one sentence, e.g., when the points of Wheeler's \lq\lq superspace" are generalized to consist of ``super Riemann manifolds'' whatever the latter are. Although there are many works published on this topic, the precise definition has to be clarified, cf. \cite{BGLS} where the importance of the non-holonomic (nonintegrable) distribution in the models of super Minkowski spaces and superstrings is demonstrated, and the local invariants of real-complex (super)manifolds --- analogs of the Nijenhuis tensor --- are defined and computed in several cases.

Possible applications of the language of supersymmetry still remain to be explored. 

Successful applications of supersymmetry in models of solid state physics (see \cite{Ef} and later works by Efetov with co-authors) are almost invisible to mass media in the shadow of discussions ``Is there supersymmetry?" concerned exclusively with high energy physics and expensive projects like the Large Hadron Collider. Whatever the discoveries of the LHC will be, the majority of researchers are certain that if Einstein's dream ``to see a~Grand Unified Theory of all interactions''  (SUSY GUT) will ever be realized, it will be formulated in the language  of supersymmetry. 

Wess and Zumino were the first to lucidly and consciously demonstrate the importance of supersymmetry, more precisely --- of Lie superalgebras, in modern terms.   Lie supergroups corresponding to these Lie superalgebras were introduced by several authors independently of my --- completely inaccessible --- preprint published in Russian in Proceedings of XIII National USSR students' scientific conference in Novosibirsk, also in 1974. Namely, a~\textit{supergroup} (not necessarily Lie) is a~ group object in the category of supermanifolds or superschemes. 

The definition of algebraic and Chevalley supergroups was given much later, but the  functorial definitions of Lie superalgebras published until recently do not work over ground fields of characteristic 2 or~3, and therefore to define algebraic and Chevalley supergroups in these characteristics is still an open problem. For a~definition of the Lie superalgebra in terms of the functor of points, see \cite{KLLS}, \cite{KLLS1}. (Another general comment to the currently published works on Chevalley supergroups --- actually, their elementary subgroups, is clear from the paper \cite{KPV}.)

For certain open problems of supersymmetry theory, see \cite{L1}, \cite{L2}, \cite{L3} and references therein.

\textbf{The Vinberg--Onishchik seminar} (names are ordered in accordance with the Cyrillic alphabet). In Moscow, Onishchik was a~co-leader of the joint Vinberg--Onishchik  seminar till the last days of his ability to move in town on his own. Among numerous active participants of this famous seminar I'll distinguish only the most prominent ones: D.~Alekseevsky, D.~Akhiezer, A.~Elashvili, V.~Kac, O.~Shwartsman, V.~Popov (later a~co-leader of the seminar), and the silent listener from the ``engineer stream'' --- V.\,M.~Sergeev --- who applied some of the knowledge acquired at the seminar to the very first descriptions of what is nowadays called ``econophysics''. Sergeev's ideas (see the translation of his book \cite{S}) were highly estimated by M.~Gell-Mann  (the Nobel prize laureate for a~description of quarks).

Being a~student majoring at the Chair of Algebra, where Onishchik was an Associate Professor, I frequented the  Vinberg--Onishchik seminar of the same Chair and the research seminar directed by A.\,A.~Kirillov (of another Chair, the one where Berezin and Gelfand worked, and whose seminars I attended only occasionally). At these seminars I delivered my results during the fall of 1973, in particular, the definitions of superscheme, supermanifold and Lie supergroup. (Oni\-shchik immediately proved that the bundle over the underlying group, the bundle whose sections define the structure sheaf of the Lie supergroup, is trivial; its fiber is $\fg_\od^*$. For a~ recent proof of this statement in the widest generality --- for group superschemes --- see \cite{MaZ}.) 

In 1971, V.\,G.~Kac has already applied his technique of $\Zee$-graded Lie algebras to define three series of what today is interpreted as ``simple Lie superalgebras of vector fields on the $0|n$-dimensional superpoint''. Kac delivered his result at the Vinberg--Onishchik seminar and at a~ National USSR (``All-Union'' as this was translated into ``Moscow English'' of those times) algebraic conference, see Kac's preface in \cite{DSB}.

After my talk, Vinberg and Onishchik, in chorus, suggested to me and V.\,G.~Kac to join forces and classify simple finite-dimensional Lie superalgebras over $\Cee$. And so we did, see acknowledgement to me ``for constructive help'' in \cite{Kac}, where the main idea of the proof and it implementation was due to Kac. For Kaplansky's accounts of his and other researchers contribution to the classification, see \cite[Newsletters]{Kapp}, \cite{FK}. (Proof of the completeness of known deformations with odd parameter existed only as a~folklore and a~conjecture until recently, see \cite{L3}.) For the classification of simple infinite-dimensional Lie superalgebras of vector fields with polynomial coefficients over $\Cee$, see still unfinished preprint of the review \cite{LSh1}, updated in \cite{BGLLS}, with main ingredient in \cite{Sh5, Sh14}. 

\textbf{Yaroslavl}. In 1975, having lost hope to be promoted to full professor at Moscow State University, Onishchik accepted an offer from Yaroslavl State University. For 30-odd years he commuted from Mos\-cow to Yaroslavl twice a~month (first, for 3-day visits; later, for a~week each time) to deliver lectures, and direct course works and Ph.D. theses. 

The distance is ca.~280 km. Today, it takes 3~hours by train (and an hour through Moscow by bus and underground); in the 1970s and 1980s, the train took a~big portion of the night. Before 2003, when Onishchik's son bought a~ flat in Yaroslavl for his father and mother to stay during their visits with comfort and decency commensurable with Onishchik's rank and title, A.L. stayed in the dorm for Ph.D. students and younger faculty (with  common kitchen, common shower, and common  toilet at the end of the corridor). 

For at least the first academic year, A.L. shared a~room in this dorm with Professor V.~Yefremovich, very interesting person (perhaps, too interesting to share a~room with, after all his tragic experiences, and with habits he acquired in detention at the time of Stalin's purges) and a~remarkable mathematician, see \url{https://en.wikipedia.org/wiki/Vadim_Yefremovich}.

\textbf{Supersymmetry}. From 1972 on, with my ``course work'' (``student research project'') as the trigger, Onishchik got more and more involved in the study of supersymmetries. After 1986, his study of Lie superalgebras and invariants of analytic supermanifolds prevailed: of 112 (a bit fewer: the list of Onishchik's published works given in the paper commemorating his 70's birthday --- \cite{AJ} --- contains redundancies (preprint=paper in English): [68J=85J], [75J=89J], [77J=88J],  [86J], for a short version, see [97J], [98J=101J], where ``J'' is for jubilee); to which I add 
published or at least preprinted works of Onishchik 
(excluding translations and books edited)  devoted to various aspects of supersymmetry theory. The bibliography given below 
continues the list of Onishchik's works published in \cite{AJ}, and in a~recent obituary \cite{AOb}. 

After 1975, almost all Onishchik's research students worked on various aspects of supersymmetry.

Onishchik's selection of the narrow topic of research in supersymmetry for himself and his students (and the choice of a~ low-profile synonym ``supercurve'' for ``superstring'' in their papers) was partly due to his 
wish to explore his own technique without elbowing his way through the competing colleagues. 

\textbf{Basic notions}. Recall that a~\textit{supermanifold} is a~ringed space, i.e., a~pair consisting of a~topological space (manifold) $M$, and a~sheaf of supercommutative rings or algebras on $M$. A~ \textit{supervariety} is a~ supermanifold with singularities. \textbf{Locally}, the \textit{sheaf $\cL_{\Lambda (\textbf{E})}$ of sections of the exterior algebra of a~given vector bundle} $\textbf{E}$ over a~manifold $M$ is completely described by the superalgebras of functions $\cF_\textbf{E}(U):=\cF\otimes \Lambda (V)$, where $\cF:=\cF(U)$ is the algebra of functions on every open domain $U\subset M$ and $V$ is the fiber of the bundle $\textbf{E}$. Morphisms of the ringed spaces 
\[
F=(\varphi, \varphi^*):\cM=(M, \cL_{\Lambda (\textbf{E})})\tto\cM'=(M', \cL_{\Lambda (\textbf{E}')}),
\]
where $\varphi:M\tto M'$ is a~morphism of manifolds, are determined by the  homomorphisms ${\varphi^*:\cF_{\textbf{E}'}(U')\tto\cF_\textbf{E}(U)}$ of the algebras of functions on the open subsets $U'=\varphi(U)\subset M'$ for any $U\subset M$. 

There are two different categories of supermanifolds: 

(1) in the category $\textsf{SMan}$, which is considered in practically all works by mathematicians and physicists since 1974, morphisms  ${\varphi^*}$ preserve parity;  

(2)  in the other category, denote it $\textsf{GMan}$, morphisms ${\varphi^*}$ are determined by arbitrary homomorphisms of the supercommutative rings of functions, see \cite{L0}, \cite{LS1}, \cite{LS2}, \cite{Lsos}. Interestingly, if $\sdim \cM= a|2b$, the integral over $\cM$ is invariant under the changes of variables that do not preserve parity. From the very beginning I conjectured that these more general symmetries are no less reasonable than those, now conventional, infinitesimally embodied by Lie superalgebras. I will not digress to this very tempting subject here, I only advise to read  \cite{Lun}, \cite{Iy}, and references therein, and my recent paper \cite{L4} where I formulate my other conjecture extending Berezin's \eqref{BerConj} from ``Grassmann" to ``Clifford", or more generally ``graded-commutative".

Clearly, there are as many objects in the category of smooth supermanifolds $\textsf{SMan}$ as there are objects in the category of vector bundles, but  in $\textsf{SMan}$ there are more morphisms than there are 
 morphisms of the vector bundles $\Lambda f: \Lambda (\textbf{E})\tto \Lambda (\textbf{E}')$, where $f: \textbf{E}\tto \textbf{E}'$ is a~ morphism of bundles, that determine these supermanifolds, see \cite[Ch.~4]{Lsos} or  \cite[Ch.~1]{Del}. For a~nice and brief introduction to the theory of ringed spaces, see~\cite{MaS}. This theory is needed to describe superschemes and supervarieties (supermanifolds with singularities) over fields of any characteristic~ $p$. 

In the categories of complex-analytic supermanifolds and superschemes, the objects  of the form $\cM:=(M, \cL_{\Lambda (\textbf{E})})$ can have deformations and the deformed objects are not of this form, so these categories  have not only more morphisms than the category of vector bundles, they have more objects. These objects are  ringed spaces $\cM:=(M, \cO_{\cM})$, where $(M, \cO_{M})$ is a~scheme or an analytic manifold, and $\cO_{\cM}$ is a~sheaf of supercommutative superalgebras or super rings. 

Let $\cI\subset  \cO=\cO_{\cM}$ be the subsheaf of ideals generated by the 
subsheaf of odd sections $ \cO_\od$ and let  
$ \cO_{\rd}:=  \cO/\cI$. 

The supermanifold 
$\cM_{rd}:=(M,  \cO_{rd})$ will be called \textit{oddly reduced} to distinguish from
the \textit{reduced} (super)manifold or (super)scheme $\cM_{red}$ whose structure sheaf is obtained by factorization of $ \cO$ modulo the subsheaf generated by \textbf{all} nilpotents (not only odd ones); the notation $\cM_{rd}$ is due to Manin, see \cite[Subsection~4.1.3]{MaG}.

Consider the following filtration of 
$ \cO$ by powers of $\cI$:
\[
 \cO=\cI^0\supset\cI^1\supset\cI^2\supset\ldots\supset\cI^n
\supset\cI^{n+1}=0.
\]
The graded sheaf ${\rm gr}\, \cO=\oplus_{0\leq i\leq n}\ {\rm gr}_i \cO$ with ${\rm gr}_i \cO:=\cI^i/\cI^{i+1}$
defines the split supermanifold $(M,{\rm gr}\, \cO)$ 
called the \textit{retract} of $(M, \cO)$.

Let $\gr \cO_{\cM}\simeq \cL_{\Lambda(\textbf{E})}$ be the associated sheaf of graded supercommutative superalgebras, see \cite{Gr},  \cite{Pa}, \cite{Va1}, \cite{MaG}. 

Any supermanifold isomorphic to ${\cM:=(M, \gr  \cO_{\cM})}$ is called \textit{split}. Observe that any smooth or analytic supermanifold is locally split. A~locally split, but not split supermanifold is called \textit{non-split}. 

The first published proof  of the statement ``all smooth supermanifolds are split'' is due to Gaw\c{e}dzki, see \cite{Ga}. This statement became known as  \textit{Batchelor's theorem}, see \cite{Bat}. It was known to participants of the Kirillov, Vinberg--Onishchik, and Manin seminars since ca 1973; for its short proof which follows from the existence of a~partition of unity in the $C^\infty$ case, see \cite[Subsection~4.1.3]{MaG}. 

\textbf{The study of non-split supermanifolds} is what Onishchik was doing with his students for the case of super Grassmannians. Lately, this topic drew attention of physicists. Witten --- one of the most prolific modern theoretical physicists and enthusiast of superstring theory --- published a~ review of the related problems and results that filled the whole issue (607~pp.) in one volume of the journal ``Pure and Applied Mathematics Quarterly", see \cite{W}.

Let me comment here on the works of Onishchik on supersymmetries, pointing at their influence  leading to generalizations made by other researchers, and indicate several open problems. The references of the form $[nJ]$, where $n\leq 101$, are repeated from the list in \cite{AJ}.

Onishchik told me that he understood the difference between the split and non-split supermanifolds having read two papers: (a) the paper  \cite{Gr} by P.~Green, who was the first to point at this difference, and lucidly describe the obstructions to non-splitness (corresponding to only even parameters)  and (b) the rich with results paper by Vaintrob \cite{Va1}. These papers encouraged Onishchik to work on non-abelian analog of Dolbeault's theorem and its applications to classification of homogeneous supermanifolds, see Onishchik's papers [93J], [96J] and \cite{104}. Let me note that Palamodov's description in
\cite{Pa} submitted only 2 month later than  \cite{Gr} and expounded in \cite[Ch.4, \S4, Sections 6--9]{Ber1} reproduced in \cite[Ch.3, Th.2, p.126]{Ber2}, as well as in \cite[Ch.4, \S2, Prop.~9, p. 191]{MaG}, are more lucid to me than Green's description which was used by Onishchik and his students. Nobody considered odd parameters of deformations of split supermanifolds into non-split ones, this gap is repaired in \cite{L3}.

$\bullet$ Observe that in all the papers on non-splitness of supermanifolds published so far, it is often assumed that the split supermanifold $\cM$ to be deformed into non-split is a~complex one. I do not see where this is needed. The underlying manifold $M$ does not have to be complex, either. We only need an \textbf{almost} complex structure of $\cM$, sometimes even less: we need a~\textbf{real-complex} structure of $\cM$, see \cite{BGLS}, where the analog of the Nijenhuis tensor in presence of a~non-integrable distribution is defined (and calculated for several types of superstrings and $N$-extended super Minkowski spaces).


[50J]=\cite{111} was deposited in VINITI (Russian Institute for Scientific and Technical Information) and is inaccessible together with all other depositions; I translated and edited it for this collection. 
In a~recent paper by D.\,V.~Alekseevsky and A.~Santi \cite{AS}, the results of  [50J] are not rediscovered, despite the similarity of the titles. 


[51J], together with a~classical result of Dynkin, was generalized in the paper \cite{ShM} by I.~Shchepochkina. Shchepochkina announced a classification of maximal simple subsuperalgebras in simple vectorial Lie superalgebras in Reports of Stockholm University no.~32/1988-15; for details, see \cite{LSh}.


[55J], [81J], [83J], and \cite{104} contain generalizations of certain parts of \cite{Va1}. In particular, for the isotropic super Grassmannians of maximal type associated with the ortho-symplectic supergroup $\mathcal{OS}p$  preserving a~non-degenerate symmetric \textbf{even} bilinear form over $\Cee$, the cohomology enabling one to determine which of these super Grassmannians is split are calculated.  An analog of this result is obtained for the isotropic super Grassmannians of maximal type associated with the \textit{periplectic}, as A.~Weil suggested to call it, supergroup $\mathcal{P}e$ preserving a~non-degenerate symmetric \textbf{odd} bilinear form over $\Cee$.

E.\,Vishnyakova \cite{Vi2} calculated
Lie superalgebras of global holomorphic vector fields for $\mathcal{GL}$- 
and   $\mathcal{GQ}$-flag supermanifolds of any type, and also for
$\mathcal{OS}p$- and $\mathcal{P}e$-flag supermanifolds of maximal type. For this, Vishnyakova had to improve results due to Onishchik and Serov \cite{104} by offering convenient atlases. \textbf{For $\mathcal{OS}p$- and $\mathcal{P}e$-flag supermanifolds of non-maximal type, the problem is still open}.


[58J] contains a~ proof of the fact that the Cartan subalgebras in any simple finite-dimensional Lie superalgebra $\fg=\fg_\ev\oplus\fg_\od$ over $\Cee$ are conjugate under the action of the group $\text{Aut}_e\fg$ generated by automorphisms of the form $\exp(\ad_x)$, where $x \in\fg_\ev$ and $\ad_x$ is nilpotent. This result rediscovers the same result in the paper by Scheunert \cite{Sch} who also described the center of $U(\fpe(n))$ (this center was also described by A.N.Sergeev at about the same time or earlier; unpublished). 

The importance of the center of $U(\fg)$ in the study of $\fg$-modules is high as is well-known. Therefore, it is worth to note that in \cite{Ser}, Serganova suggested to replace $U(\fspe(n))$, whose center is trivial, by $\overline U(\fg):=U(\fg)/\fr(U(\fg))$ whose center is sufficiently big, where $\fr(A)$ is the radical of the algebra $A$. It is very tempting to apply Serganova's idea to other Lie (super)algebras $\fg$ with trivial center of $U(\fg)$ (and with non-trivial centers as well).

In  the 1970s,  several researchers defined a~Cartan subalgebra of a~given Lie superalgebra $\fg$  as any Cartan subalgebra of $\fg_\ev$. Under this definition the statement on conjugacy of Cartan subalgebras is trivially equivalent to the same statement on Cartan subalgebras of Lie algebras (provided the ground field is algebraically closed of characteristic 0) and does not reflect the specifics of ``super''. 

In \cite{Sch} and [58J] a~\textit{Cartan subalgebra} $\fh$ of a~given Lie superalgebra $\fg$ is defined as a~nilpotent Lie superalgebra coinciding with its normalizer in $\fg$. In \cite{Sch} and [58J], the authors show that the Cartan subalgebras in (finite-dimensional) $\fg$ and $\fg_\ev$ are in one-to-one correspondence: for any Cartan subalgebra $\fh$ of $\fg$, its even part $\fh_\ev$ is a~Cartan subalgebras in $\fg_\ev$, and $\fh$ coincides with the weight space of weight 0 in the adjoint representation of $\fh_\ev$ in $\fg$. The proof is based on the following analog of Engel's theorem:

\textit{Let $\rho: \fg\tto\fgl(V)$ be a~linear representation, $V\neq 0$ and all operators $\rho(X)$, where $X\in\fg_\ev$ are nilpotent. Then, there exists a~non-zero $v_0\in V$ such that $\rho(X)v_0=0$ for all $X\in\fg$.} 

In turn, the proof of this analog of Engel's theorem is based on the fact ``$\fg$ is solvable if and only if $\fg_0$ is solvable''; for its proof over $\Cee$, see \cite{Sg}. Observe that  if $\Char\Kee =2$, there are examples of simple Lie superalgebras $\fg$ with solvable $\fg_\ev$, see \cite{BGL}, \cite{BGLLS}. 

There are many cases where the Lie superalgebras (even finite-dimensional over $\Cee$) behave like Lie algebras over fields of positive characteristic or as infinite-dimensional Lie algebras. Therefore, the result of \cite{Sch} and [58J] on conjugacy of Cartan subalgebras is not expected.
In \cite{Pre},
Premet  proved that if $\Char\Kee>5$, then the Cartan subalgebras of simple Lie algebra $L$ operate on $L$ via upper triangular matrices and gave a~counterexample for $\Char\Kee=5$. These facts had been useful in the classification of simple finite-dimensional modular Lie algebras over algebraically closed fields $\Kee$ such that $\Char\Kee>3$; no doubt they will be used in the future classification of simple modular Lie superalgebras.

For a~latest discussion of what should we take for a~definition of the Cartan subalgebra (with some answers), see \cite{Ray}.

I.~Penkov with co-authors investigated what is the ``right'' definition of the Borel subalgebra of a~(simple, or close to it, like $\fgl$, $\fpgl$ and $\fsl$  are close to $\fpsl$) Lie superalgebra; for a summary of their studies in infinite-dimensional cases, see \cite{HP}. There are two natural candidates: 

$\bullet$ A~maximal solvable subalgebra (as classified in \cite{Shch} for the $\fgl$ and $\fsl$ series).

$\bullet$ A~smaller subalgebra, namely, the semi-direct sum of a~Cartan subalgebra and the nilpotent subalgebra spanned by positive root vectors, see \cite{DPS}, \cite{DP}, \cite{PS}, \cite{GY}, \cite{HZ}. This subalgebra seems more useful from the  point of view of representation theory than the maximal one. Observe here that the notion of roots of Lie superalgebras should be defined not as for Lie algebras over fields of characteristic 0; for a~correct definition suitable for Lie superalgebras and modular Lie algebras, see \cite{BLLoS}.

Since \textbf{in various aspects, finite-dimensional Lie superalgebras over $\Cee$ resemble modular, and infinite-dimensional Lie algebras, and quantum algebras (``groups")} (for example, see \cite{J}), it is interesting to apply these results by Penkov and his co-authors to the modular Lie (super)algebras, and to the quantum versions. 

[68J=85J] Here, V.\,A.~Bunegina and A.\,L.~Onishchik classified, up to an isomorphism, homogeneous complex supermanifolds of dimension $1|m$ with the underlying manifold $\Cee\Pee^1$   for the case where $m \leq 3$. 

To formulate their result, consider a~ cover of
${\mathbb{CP}}^1$ by two charts $U_0$ and  $U_1$ with local coordinates $x$ and  $y={x}^{-1}$, respectively. Then, for
${\mathcal{CP}}^{1|3}_{k_1k_2k_3}=(\Cee\Pee^1, \cL_{\Lambda (\textbf{E})})$, 
the transition functions 
in $U_0\cap U_1$ are 
$y ={x}^{-1}$ and $\eta_i ={x}^{-k_i}\xi_i$ for $i=1,\dots,3$ and for $k_i\in\Zee_{\geq 0}$,
where the $\xi_i$ and $\eta_i$ are basis sections of the bundle $\textbf{E}$ over
$U_0$ and $U_1$, respectively.

For $m = 1$, many authors claimed that any supermanifold of odd dimension 1 is split. This is not true, actually: take into account odd parameters of deformations, see \cite{L3}.

For $m = 2$, Bunegina and  Onishchik showed that there exists only one non-split homogeneous supermanifold, the one constructed by P.~Green (see \cite{Gr}). Observe that Manin's answer \cite[Ch.4, \S2, Prop.~9, p. 191]{MaG} is at variance with that by Bunegina and Onishchik and I intend to answer in a separate paper whose answer is correct and where is the mistake.

For $m = 3$, Bunegina and  Onishchik showed that there exists a~series of non-split homogeneous supermanifolds, whose retract correspond to the sum of line bundles of the form $L_{-k_1}\oplus L_{-k_2}\oplus L_{-k_3}$, where $k_1=k_2=2$ and either $k_3=0$ or $k_3\geq 2$. For any such triple  $(2,2,k_3)$, there exists only one (class of)  non-split homogeneous supermanifolds.

Recall that the definition of the action (in particular, a~transitive one), of the Lie supergroup $\cG$ on the supermanifold $\cM$ is a~ natural superization of the Lie group action from algebraic geometry. The supermanifold $\cM$ is called \textit{even-homogeneous} (\textit{$\ev$-homogeneous} for brevity) if $\cG$ is a~group, not a~supergroup. In \cite{Vi1}, E.\,G.~Vishnyakova generalized the results of [68J=85J], together with the pioneer P.~Green's result (see \cite{Gr}), by describing non-split $\ev$-homogeneous supermanifolds of superdimension $1|3$ whose underlying manifold is the projective line. 


[77J=88J] Here, the paper \cite{FN}  is interpreted from the point of view of supermanifolds and the results are used to calculate the cohomology of the projective space $\Cee\Pee^n$ with coefficients in the sheaf of vector-valued differential forms. Grozman  (see~\cite{Gz}) also interpreted the space of differential forms on a~manifold $M$ with values in the Lie algebra of vector fields on $M$ as the Lie superalgebra which is the centralizer of the exterior derivation operator $d$ considered as a~vector field on the supermanifold $(M, \Omega)$ whose structure sheaf is the sheaf of differential forms on~ $M$.


[90J, 91J] Onishchik and Platonova continued the study of homogeneous and $\ev$-ho\-mo\-ge\-ne\-ous complex supermanifolds of superdimension $n|m$, where $n\geq 2$ and $m\leq n$, and the underlying manifold $\Cee\Pee^n$: all non-split supermanifolds of this type, whose group of automorphisms is of the $PGL$ type, are classified.   (In the non-split case, it is assumed that the retract of $(M, \cO)$ is $\ev$-homogeneous relative to $\mathfrak{pgl}(n+1)\simeq\fsl(n+1)$.)


[92J, 93J, 100J], \cite{107}, \cite{108} These are reviews suitable for an initial acquaintance with the problems of description and classification of non-split supermanifolds; the paper [108J] contains, moreover, new classification results.


\cite{109} N.~Ivanova and Onishchik described parabolic subalgebras of simple or ``close'' to simple (such as $\fgl$ or $\fq$ or $\fpe$ are close to $\fsl$ or $\fpsl$, or $\fpsq$, or $\fspe$, respectively) Lie superalgebras where the \textit{parabolic } subalgebra of any (finite-dimen\-si\-o\-nal) $\Zee$-graded Lie superalgebra $\fg$ is defined as $\fp:=\oplus_{i\geq 0} \fg_i$. 
These definitions and results are used in the papers \cite{Srg}, \cite{GY}. 

Observe that \textbf{the supermanifolds $\cC\cP^{1|n}$, and their deformations, are known among physicists as superstrings}. I am sure that if this observation had been clearly written in each of the papers by Onishchik and his students, 
the interest to the results  contained in these papers would have been much more vivid\dots

The purpose of the paper \cite{108} is to classify (up to an isomorphism) $\ev$-homogeneous non-split complex supermanifolds of dimension $1|m$, where $m\leq 3$, whose underlying manifold is $\Cee \Pee^1$. The answers are as follows.

For $m = 1$, these supermanifolds are split. In \cite[Ch.4, \S2, Prop.~8, p. 190]{MaG} and \cite[Example 3.3.1 (1)]{Va1}, there is formulated that, moreover, all supermanifolds $\cM=(M, \cO)$ of superdimension $m|1$ are split. However, this is true only if we consider $\cM=(M, \cO)$ ``individually", not over a~supervariety of parameters, as other deformation problems were considered in \cite{Va1} and \cite{MaG}. For examples of non-zero obstructions, see \cite{L3}. 

For $m = 2$, Bunegina and Onishchik claimed that there exists only one such supermanifold constructed in 1982 by P.~Green (see \cite{Gr}): it is the result of a~deformation of 
\[
\text{$\cC\cP^{1|2}:=(\Cee\Pee^1, \cL_{\Lambda(\textbf{E})})$, where $\textbf{E} = 2T^*(\Cee\Pee^1):=T^*(\Cee\Pee^1)\oplus T^*(\Cee\Pee^1)$. }
\]
As I mentioned above, this is at variance with Manin's answer \cite[Ch.4, \S2, Prop.~9, p. 191]{MaG}.

For $m = 3$, there exists a~series of non-split $\ev$-homogeneous supermanifolds, parameterized by elements in $\Zee \times \Zee$, three series of non-split $\ev$-homogeneous supermanifolds, parameterized by elements of $\Zee$, and a~finite set of exceptional supermanifolds. 

Observe that in all the works listed in the bibliography, the Lie superalgebras and the classification of non-split supermanifolds are considered ``naively'', i.e., in terms of geometric points. 

 \textbf{Open Problem: Consider Lie superalgebras and the classification of non-split supermanifolds functorially, i.e., with odd parameters}. This problem was also formulated in \cite[Ch.3, p.136]{Ber2}, \cite{Va1}. For examples of its solution, see \cite{L3}, \cite{L2}. 

\textbf{Research students}. In the province, with very few students interested in research, Onishchik managed, nevertheless, to encourage several of them, E.~Vishnyakova being the most active and successful. The others did not manage, regrettably, withstand the ``resistance and pressure of the environment'' (more precisely,  find a~ way out of poverty while continuing pure research in supersymmetry when salaries of faculty members at the Universities became lower than those of janitors and either of these salaries became insufficient for survival, to say nothing about decent life).  As one can judge by the lack of publications, they stopped their research or switched from the study of supersymmetries to ostensibly more applied problems.

From the total number of about 30 of Onishchik's students who have defended ``kandidatskaya'' (Ph.D.) theses, among whom several have defended the second (Dr. Sci. a.k.a. Habilitation) thesis, I'll list the students from Yaroslavl and nearby towns Tver and Rybinsk who actively studied supersymmetry. Mostly, the works of Onishchik's students during the past 30 years were devoted to solution of various aspects of the following problem explicitly formulated in [82J]:  
\be\label{pr}
\begin{minipage}[l]{11.5cm}
``for a~given split supermanifold $(M, \cO)$, classify, up to an isomorphism, all homogeneous, or at least $\ev$-homogeneous, 
complex supermanifolds whose retract is $(M, \cO)$ with a~given homogeneous 
complex manifold $M$''.
\end{minipage}
\ee

In [82J], Onishchik classified all supermanifolds with retract $(M, \Omega)$, where $M$ is an irreducible simply connected compact Hermitian symmetric space, $\Omega$ is the sheaf of holomorphic forms on $M$. In particular, Onishchik proved that the only homogeneous such supermanifold is the $\Pi$-symmetric super Grassmannian (first described in \cite{MaG}).

In [74J], Bunegina and Onishchik  constructed a~1-parameter family of homogeneous supermanifolds whose retract is the superstring (or supercurve) $\cC\cP^{1|4}$. The Lie superalgebra of vector fields on this supermanifold is a~ representative $\Gamma(\sigma_1,\sigma_2,\sigma_3)$  of a~parametric family  $\fosp_\alpha(4|2)$  of deformations of $\fosp(4|2)$ defined for the triples ${\sigma\neq (0,0,0)}$ such that 
\[
\sigma_1+\sigma_2+\sigma_3=0. 
\]
This family was discovered by Kaplansky, see \cite{Kapp}; now it is often called $D(2,1;\alpha)$, where assuming (up to renumbering) $\sigma_1=1$ we set $\alpha=\frac{\sigma_i}{\sigma_j}$ for $\sigma_j\neq 0$, where $\{i,j\}=\{2,3\}$. Unlike $\Gamma(\sigma_1,\sigma_2,\sigma_3)$  or $\fosp_\alpha(4|2)$, the name 
$D(2,1;\alpha)$ is ill-chosen as explained, for example, in \cite{CCLL}.

The classification of superstrings with retract $\cC\cP^{1|m}:=(\Cee\Pee^{1}, \cO_{\Cee\Pee^{1}}\otimes E^{\bcdot}(\textbf{E}))$ for $m < 6$ reduces to the calculation of \v{C}ech 1-cohomology of $\Cee\Pee^{1}$ with coefficients in  the tangent sheaf of $\Cee\Pee^{1}$ tensored by an exterior power of the bundle $\textbf{E}$, see \cite[Ch.3]{Ber2}, \cite{Va1}.

E.~Vishnyakova and M.~Bashkin constructed explicit bases of these cohomology spaces and  bases of $\fsl_2(\Cee)$-invariant cohomology classes.  

If $m = 3$, Vishnyakova gave classification (up to an isomorphism) of $\ev$-homogeneous superstrings with the underlying manifold  $M=\Cee\Pee^1$, see \cite{Vi1}.

If $m = 4$, Bashkin gave classification (up to an isomorphism) of $\ev$-homogeneous superstrings with the underlying manifold 
$M=\Cee\Pee^1$, see \cite{Ba1}, \cite{Ba2}. 

If $m = 5$, Bashkin described $\ev$-homogeneous superstrings with retract $\cC\cP^{1|5}$ and found systems of equations for the coefficients of linear combinations of basis 1-cocycles. These systems (\textbf{no solution of any of them is found yet})  is contained in Bashkin's Ph.D. thesis; they single out homogeneous superstrings with this retract, see \cite{Ba5}. 

In the case where the retract $\cT^{m|n}$ of $\cM$
is associated with the trivial bundle of rank $n$ over the $m$-dimensional complex torus $M=T^m$, Bashkin proved that $\cM$ is always $\ev$-homogeneous; it is homogeneous if and only if $\cM\simeq\cT^{m|n}$, see \cite{Ba3}.

Onishchik studied a~similar classification problem 
for the case where the underlying manifold is ${M=Gr_2^4}$, and formulated the answer if the representation of the stabilizer of the point defining the retract is completely reducible. This was later generalized by Igonin, see \cite{Ig}, \cite{Ig2}. Igonin's results, though reviewed in English in \cite{106}, \cite{107}, deserve more accurate study and ``rounding up''. 

Solving the same problem, Serov calculated 1-cohomology with coefficients in the tangent sheaf for several split homogeneous supermanifolds  whose underlying manifolds are flag manifolds.

Onishchik and Vishnyakova calculated the Lie superalgebras of vector fields for the 4 series of flag supermanifolds related with Lie superalgebras of series $\fgl$, $\fq$, $\fosp$, $\fpe$. For example, the main results of \cite{109} are the following: 

$\bullet$ The
classification of locally free sheaves of modules over supermanifolds with a~given
retract in terms of non-abelian 1-cohomology. 

$\bullet$ The study of locally free 
sheaves on projective superspaces, in particular, a~superization of
the Barth--Van~de~Ven--Tyurin theorem (any finite rank vector bundle on $\Cee\Pee^\infty$ is isomorphic to a~ direct sum of line bundles). 

$\bullet$ A spectral
sequence which helps to calculate  cohomology with values in a~locally free
sheaf of modules. In
the case of the split supermanifold, the necessary and sufficient
conditions for triviality of cohomology class which corresponds to
the tangent sheaf are given.

The paper \cite{Vi0} is clearly written and is often cited.
Results of other papers by Vishnyakova:

$\bullet$  Proof of rigidity of flag supermanifolds under certain rigidity assumptions on the flag supermanifold, i.e., its complex structure does not admit any non-trivial small deformation. (Moreover, under the same assumptions, the~flag supermanifold is proven to be a~unique non-split supermanifold with a~given retract.)

$\bullet$ Vector bundles and \textit{double vector bundles}, a.k.a. $2$-fold vector bundles, arise naturally, for instance, as base spaces for algebraic structures such as Lie algebroids, Courant algebroids and double Lie algebroids. (Recall in this instance that Vaintrob showed that all these structures can be uniformly descried in the language of super\-geometry, see \cite{Va}, in terms of the \textit{homological} vector field, i.e., an odd vector field $D$ such that $D^2=\frac12[D, D]=0$.) Vishnyakova established a~link between the super and classical pictures by the geometrization process, leading to an equivalence of the category of graded manifolds of degree $\leq 2$ and the category of (double) vector bundles with additional structures. 

$\bullet$ In \cite{Vi3}, Vishnyakova established an equivalence between a~subcategory of the category of $n$-fold vector bundles and the category of graded manifolds of lattice type $\Delta$ (for definitions of the non-common notions involved, see \cite{Vi3}).

A. V. Sudarkin studied gradings and parabolic subalgebras of Lie superalgebras ``of classical type'' in terms of root systems of  simple finite-dimensional Lie algebras; for the answer in the case of $\fsl(m|n)$, see \cite{Su}. 

A.\,A.~Serov (see \cite{Sr})  described geometric points of spinor supergroups which are two-sheeted covers of the ortho-symplectic supergroup. For an exposition of his results with suggestions as to how his result should be generalized in several directions, see \cite[pp.~397--401]{Lsos}, a~bit more accessible than the Yaroslavl collection containing \cite{Sr}. (Although the English version of this  result of Serov had not been published yet, the paper  \cite{Sr} written in Russian was reviewed, and the review is unusually lucid, see MR0881845. Nevertheless,  a~less clear description of $\Ree$-points of spinor supergroups was published 30 years later, see \cite{SAS}.)   

These are ``direct" superizations of spinor groups.
Observe that there is also an analog of the spinor group which is a~ two-sheeted cover of the periplectic supergroup, whose $\cC$-points are described in \cite[Ch.1]{Lsos}.

Serov's other ``super'' results concern deformations of complex supermanifolds and description of vector fields on flag supermanifolds of maximal type, see \cite{104}. 
 
 O.\,M.~ Sulim's results on maximal solvable Lie subsuperalgebras of $\fgl(m|n)$ and $\fsl(m|n)$ were completed by I.\,M.~Shchepochkina who gave a~ classification, see  \cite{Shch}.  
 
G.\,V.~Egorov's results \cite{E},  as well as partial results of other researchers on normal forms (shapes) of multivectors, were interpreted as results on normal shapes of germs of  functions on the superpoint and generalized by Sergano\-va and Vaintrob in the paper \cite{SV}. This paper contains the description and classification of germs of singularities up to the action of the group of diffeomorphisms on the superdomain. To classify germs of singularities up to the action of the \textbf{super}group of diffeomorphisms remains an \textbf{open problem}. A version of the problem Egorov solved is considered in \cite{BdH}.

In \cite{E1}, Egorov showed that there are several types of  Lie superalgebras of infinite supermatrices in each of the types $\fsl$, $\fosp$, $\fpe$ and $\fq$, some of them having a~non-trivial central extension, some having a~supertrace (important properties, especially from the point of view of mathematical physicists). During the past 4 decades Penkov with co-authors studied various types of infinite supermatrices forming Lie superalgebras of the series $\fsl$, $\fosp$, $\fpe$ and $\fq$, but not of the types distinguished by Egorov. These Lie superalgebras remain to be studied.

A.\,V. Yastrebov's paper \cite{Ya} on the Cramer  theorem and the  Cayley-Hamilton theorem for matrices over a~superalgebra was generalized in several directions, see a~review \cite{LL} and references therein.


\footnotesize
\textbf{Editing}. In addition to editing translations mentioned in the Fest\-schrift \cite{AJ}, Onishchik edited other books, here is as much as I could find, but certainly not the complete list: 

M. Noeter ``Sophus Lie''. (Russian) Translated from the German by B. R.~Frenkin. Edited by A. L.~Onishchik. Istor.-Mat. Issled. (2) No. 11(46) (2006), {306--347}, 359. 

 Dynkin E. B. ``Selected papers of E. B. Dynkin with commentary''. Edited by A.\,A.~ Yushkevich, G. M. Seitz and A. L. Onishchik. American Mathematical Society, Providence, RI; International Press, Cambridge, MA, 2000, xxviii+796~pp. 
 
Grauert H., Remmert R. ``Theorie der Steinschen R\"aume''. (German) [Theory of Stein spaces] Grundlehren der Mathematischen Wissenschaften [Fundamental Principles of Mathematical Sciences], 227. Springer-Verlag, Berlin-New York, 1977. xx+249~pp. Translated into Russian by D. Akhiezer, translation edited and appended by A.~Onishchik, Nauka, Moscow, 1989, 336~pp.;
 
Grauert G., Remmert R. ``Analytic local algebras''. With the collaboration of O. Rimenshneider. Translated from the German by D. N. Akhiezer. Translation edited by A. L. ~Onishchik. "Nauka'', Moscow, 1988, 304~pp.

Onishchik invested (I can not write ``spent'', although it was mainly others who harvested his investments) enormous amounts of time into aperiodic collections of scientific works which he not only edited completely alone, but also often checked galleys instead of the authors. Partly, these collections were refereed by Math Reviews making the results published there at least theoretically available world-wide (regrettably, the   review in MR is, not seldom, less informative than the Abstract written by the authors). I encourage the reader to appreciate the amount of Onishchik's work by looking only at the precious little which \textbf{has} appeared in Math~Reviews:

``Questions of group theory and homological algebra'' A.L.Onishchik (ed.), Yaro\-slavl, Yaroslavl State University, 1998. 287~pp.; 1994. 156~pp.; 1992. 168~pp.; 1991. 172~pp.; 1990. 171~pp.; 1989. 172~pp.; 1988. 172~pp.; 1987. 172~pp.; 1985. 167~pp.; 1983. 140~pp.; 1982. 164~pp.
 
``Mathematics in Yaroslavl University. Towards 20-th anniversary of the Mathematical Department''. V.\,G.~Dur\-nev, L.\,S.~Kazarin, A.\,L.~Onishcik (eds.) Yaroslavl, Yaroslavl State University, 1996. 203~pp. 
 
``Proceedings of the 5th Kolmogorov lectures''. Yaroslavl, Yaroslavl Pedagogical University, ~2007. 

Among other not covered by Math Reviews collections Onishchik (co-)edited, are ``Geometric methods in questions of algebra and analysis'', ``Problems of natural and humanitarian studies. Mathematics, informatics''. I was unable to find the details. 

\textbf{Reviewing}. Onishchik wrote a huge number of referee reports for Math Reviews and its Russian analog (RZh Mathematica). He worked on these reports almost to the very end of his life.

\normalsize

\textbf{Charity}. 
Several times Onishchik helped my friends and/or students. (I cannot write ``former students", although the rules of English allow. As I.\,M.~Gelfand used to repeat, this is like saying ``my former daughter''.) For example, Onishchik  served as official scientific advisor of those who were not recommended to the Ph.D. school by Mekh-Mat's communist party bosses, and whose actual advisor was somebody else, somebody without official license to perform such a~highly responsible and politically sensitive task as education of future Ph.D. holders.

In this way, B.\,L.~Feigin, V.\,V.Serganova, and Yu.\,Yu.~Kotchetkov defended their Ph.D. theses with Onishchik as their official supervisor. 

E.~Poletaeva almost completed all the requirements in the Ph.D. school under Onishchik's official guidance. Later, while at Penn State University, she published and defended an expounded version of what she had done in Yaroslavl. For a~summary of her results at that time, partly published in difficult to access proceedings, see \cite{Po}.

F.~Weinstein reminded me that he also was for a~while a~Ph.D. student in Yaroslavl un\-der Onishchik's official guidance but actually under D.B.~Fuchs; this helped him to get the results of compulsory tests acknowledged by the Ph.D. school  in the Netherlands when he defended his Ph.D. thesis there. 

M.~Borovoi also defended his Ph.D. thesis with Onishchik's as his official advisor after his actual initial scientific advisor --- I.\,I.~Pyatetsky-Shapiro --- emigrated. 

Later on, V.\,V.~Serganova, E.~Poletaeva, E.~Vishnyakova and Yu.\,Yu.~Kotchetkov (the first one to this day and the latter one for a~while) worked studying various aspects of Lie superalgebras. 

For her results in the theory of Lie superalgebras Vera Serganova was elected in 2018 to the American Academy of Arts and Sciences; she was twice invited to deliver a~lecture at the International Congress of Mathematicians (one of the talks was plenary).

Feigin's invited talk at ICM expounded his Yaroslavl Ph.D. thesis. 

During the post-Soviet time, when grants started to replace salaries, Onishchik incorporated numerous colleagues and students in his grants, even if the topic of their research was far from his own interests. 

Onishchik never expected any gratitude for such time-consuming help a~tiny part of which I described above.


\textbf{A story}. Onishchik had a~keen sense of humor; it helped him, and those close to him, to tolerate the reality. For example, when he was still  on the faculty of Mekh-Mat, and I was officially his student, an incident happened. Like other popular professors, Onishchik was always supervising several research students simultaneously; at that time, one of his older students was in the Ph.D. school. This Ph.D. student was in love with one of his peers; she was even his wife. However, she used to show signs of affection for other young men. 

Shortly before the time of this story, I watched one Polish movie in which the main heroine, in a~similar situation,  told one of the main male characters: ``Why do you worry?! Yes, I love him, but I love you too. Yes, I sleep with him, but I do sleep with you too!'' 

The student (like the hero in the movie) did not, however, want to share. He took a~small ax, went and hacked his wife's other partner  who lived in the same dorm. 

A~frail wonk, he hacked unsatisfactorily, as soon became clear, but being unaware of this, or nevertheless, he got very distressed and immediately threw himself out of the window. From the 21st floor. 

This nerd completely forgot in the thick of fray that 2 stories below there was a~wide roof of a~technical 19th floor on which he landed and broke his leg. 

As a~result, the femme fatale --- the cause of the trouble --- had to bring fruits and whatever to both sufferers hosted in the same hospital. The whole department and the dorm were roaring with laughter, savoring the details, except for Arkady Onishchik: he was officially reprimanded for ``failure of his pedagogical work''. This black spot in his dossier effectively annihilated his chances for promotion. The recommendation letter from such an unreliable pedagog could be considered negatively the more positive the letter was, so this love story had grim consequences for several people not directly involved at~all.


\subsection*{Acknowledgements} I was supported by the grant AD 065 NYUAD. 

\end{paper}